\documentclass[12pt,noamsfonts]{amsart}

\usepackage{amsfonts}

\def\AA{{\mathbb A}}
\def\ZZ{{\mathbb Z}}
\def\NN{{\mathbb N}}
\def\RR{{\mathbb R}}
\def\CC{{\mathbb C}}

\def\fp{{\mathfrak p}}
\def\fg{{\mathfrak g}}

\def\cG{{\mathcal G}}

\def\cU{{\mathcal U}}

\def\aa{{\bf a}}
\def\bb{{\bf b}}
\def\cc{{\bf c}}
\def\dd{{\bf d}}
\def\ee{{\bf e}}
\def\uu{{\bf u}}
\def\vv{{\bf v}}
\def\00{{\mathbf 0}}
\def\11{{\mathbf 1}}
\def\22{{\mathbf 2}}

\def\tR{{\widetilde R}}
\def\tI{{\widetilde I}}
\def\tM{{\widetilde M}}
\def\tG{{\widetilde \cG}}

\newcommand{\gr}{\operatorname{gr}}
\newcommand{\Ann}{\operatorname{Ann}}
\newcommand{\Spec}{\operatorname{Spec}}

\newcommand{\GKdim}{\operatorname{GKdim}}
\newcommand{\Kdim}{\operatorname{Kdim}}
\newcommand{\Ch}{\operatorname{Ch}}
\newcommand{\IN}{\operatorname{in}}
\newcommand{\PR}{\operatorname{PR}}
\newcommand{\GR}{\operatorname{GR}}
\newcommand{\GF}{\operatorname{GF}}
\newcommand{\Ker}{\operatorname{Ker}}

\newtheorem{lemma}{Lemma}[section]
\newtheorem{theorem}[lemma]{Theorem}
\newtheorem{corollary}[lemma]{Corollary}
\newtheorem{proposition}[lemma]{Proposition}

\theoremstyle{definition}

\theoremstyle{remark}
\newtheorem{example}[lemma]{Example}
\newtheorem{remark}[lemma]{Remark}

\begin{document}

\title[Irreducible Components of Characteristic Varieties]
{Irreducible Components of \\ Characteristic Varieties}
\author[G.~G.~Smith]{Gregory G.~Smith}
\address{Department of Mathematics, University of California, Berkeley,
California, 94720}
\email{{\tt ggsmith@math.berkeley.edu}}

\begin{abstract}
We give a dimension bound on the irreducible components of the
characteristic variety of a system of linear partial differential
equations defined from a suitable filtration of the Weyl algebra
$A_{n}(k)$.  This generalizes an important consequence of the fact
that a characteristic variety defined from the order filtration is
involutive.

More explicitly, we consider a filtration of $A_{n}(k)$ induced by any
vector $(\uu,\vv) \in \ZZ^{n}\times\ZZ^{n}$ such that the associated
graded algebra is the commutative polynomial ring in $2n$
indeterminates.  The order filtration is the special case $(\uu,\vv) =
(\00,\11)$.  Any finitely generated left $A_{n}(k)$-module $M$ has a
good filtration with respect to $(\uu,\vv)$ and this gives rise to a
characteristic variety $\Ch_{(\uu,\vv)}(M)$ which depends only on
$(\uu,\vv)$ and $M$.  When $(\uu,\vv) = (\00,\11)$, the characteristic
variety is involutive and this implies that its irreducible components
have dimension at least $n$.  In general, the characteristic variety
may fail to be involutive, but we are still able to prove that each
irreducible component of $\Ch_{(\uu,\vv)}(M)$ has dimension at least
$n$.
\end{abstract}

\maketitle

\section{Introduction}

The geometry of the characteristic variety plays a central role in the
study of systems of linear partial differential equations.  In
algebraic analysis, the characteristic variety of a linear system is
obtained from a filtration of the corresponding ${\mathcal D}$-module.
When ${\mathcal D}$ is equipped with the order filtration, Sato,
Kashiwara and Kawai~\cite{SKK} and Gabber~\cite{G} show that the
characteristic variety of any ${\mathcal D}$-module is involutive with
respect to the natural symplectic structure on the cotangent bundle.
As a consequence, they deduce the ``Strong Fundamental Theorem of
Algebraic Analysis'', which says that, if ${\mathcal D}$ is the sheaf of
differential operators on an $n$-dimensional variety, then each
irreducible component of the characteristic variety must have
dimension at least that $n$.  In this paper, we work over a field $k$
and consider left $A_{n}(k)$-modules which correspond differential
operators on $\AA_{k}^{n}$ with polynomial coefficients.  Our goal is
to extend this dimension bound to all filtrations of the Weyl algebra
$A_{n}(k)$ for which the associated graded ring is a commutative
polynomial ring in $2n$ variables and to extend this assertion to a
larger class of algebras.

We are primarily interested in this larger class of filtrations
because of its connection with monomial ideals in a commutative
polynomial ring.  Specifically, for a generic vectors $(\uu,\vv)$, the
characteristic variety $\Ch_{(\uu,\vv)}(M)$ is given by a square-free
monomial ideal or equivalently a simplicial complex.  Monomials ideals
form an important link between algebraic geometry, combinatorics and
commutative algebra.  Much of the success of Gr\"{o}bner bases theory
comes from an understanding of monomial ideals.  We believe that
further exploration of this connection will lead to new insights into
$A_{n}(k)$-modules and, perhaps, primitive ideals.  Problems of making
effective computations in algebraic analysis, provide a secondary
motivation for considering filtrations other than the standard or
order filtration.  The choice of filtration can significantly effect
the complexity of the characteristic variety.

To state our theorems more explicitly, we introduce some notation.  We
write $x_{1},\ldots,x_{n},y_{1},\ldots,y_{n}$ for the generators of
$A_{n}(k)$ satisfying the relations $x_{i}x_{j}-x_{j}x_{i}=0$,
$y_{i}y_{j}-y_{j}y_{i}=0$ and $y_{i}x_{j}-x_{j}y_{i}=\delta_{i,j}$,
where $\delta_{i,j}$ is the Kronecker symbol.  Each vector $(\uu,\vv)
\in \ZZ^{n}\times\ZZ^{n}$ induces an increasing filtration on
$A_{n}(k)$ by setting $\deg(x_{i}) = u_{i}$ and $\deg(y_{i}) = v_{i}$.
We shall focus those vectors $(\uu,\vv)$ for which the associated
graded ring $\gr_{(\uu,\vv)}(A_{n}(k))$ is the commutative polynomial
ring $S = k[\bar{x}_{1}, \ldots, \bar{x}_{n}, \bar{y}_{1}, \ldots,
\bar{y}_{n}]$.  Analogously, for any finitely generated left
$A_{n}(k)$-module $M$, we can filter $M$ by assigning degrees to a
generating set of $M$.  The associated graded module $\gr(M)$ is then
a module over the polynomial ring $S$ and, hence, the prime radical of
the annihilator of $\gr(M)$ defines a variety in $\AA_{k}^{2n}$.  This
variety is called the characteristic variety $\Ch_{(\uu,\vv)}(M)$ of
$M$.  It is independent of the choice of degrees and generators for
$M$, but depends on the filtration of $A_{n}(k)$.  The main result of
this paper is the following:
\begin{theorem} \label{main}
Let $k$ be a field of characteristic zero and let $M$ be a finitely
generated left $A_{n}(k)$-module.  If the integer vector $(\uu,\vv)$
induces a filtration satisfying $\gr_{(\uu,\vv)}(A_{n}(k)) = S$, then
every irreducible component of $\Ch_{(\uu,\vv)}(M)$ has dimension at
least $n$.
\end{theorem}
\noindent The conclusion of this theorem is vacuously satisfied when
$\Ch_{(\uu,\vv)}(M)$ is empty and this can occur as the following
example illustrates: if $M =A_{2}(k) / A_{2}(k) \cdot I$, where $I$ is
the left ideal $\langle y_{1}-1, y_{2}-1 \rangle$, then
$\Ch_{(\22,-\11)}(M)$ corresponds to the $S$-ideal $\langle 1 \rangle$
indicating that the characteristic variety is empty.  However, for
$(\uu,\vv)$ nonnegative and $M \neq 0$, the characteristic variety
$\Ch_{(\uu,\vv)}(M)$ is never empty.

Theorem~\ref{main} refines Bernstein's inequality \cite{Bo} which says
that there exists an irreducible component of $\Ch_{(\11,\11)}(M)$ of
dimension at least $n$.  On the other hand, for the order filtration
$(\uu,\vv)=(\00,\11)$, Theorem~\ref{main} follows from the fact that
the characteristic variety $\Ch_{(\00,\11)}(M)$ is involutive with
respect to the natural symplectic structure on $\AA_{k}^{2n}$.  The
involutivity of $\Ch_{(\00,\11)}(M)$ was first established by Sato,
Kashiwara and Kawai~\cite{SKK} using mirco-local analysis;
Gabber~\cite{G} provided a purely algebraic proof.  In our more
general case, a different proof is necessary because the
characteristic variety $\Ch_{(\uu,\vv)}(M)$ is not always involutive
under the natural symplectic structure on $\AA_{k}^{2n}$.  For
example, the characteristic variety of $A_{2}(k) / A_{2}(k) \cdot I$,
where $I$ is the left ideal $\langle y_{1}^{2}-y_{2},
x_{1}y_{1}+2x_{2}y_{2} \rangle$, with respect to the vector
$(1,1,1,3)$ is given by the non-involutive $S$-ideal $\langle
\bar{x}_{2}, \bar{y}_{2} \rangle \cap \langle \bar{y}_{1}, \bar{y}_{2}
\rangle$.  

The general techniques used in the proof of Theorem~\ref{main} apply
to a larger collection of $k$-algebras.  We develop these methods for
a skew polynomial ring $R$ which is an almost centralizing extension
of a commutative polynomial ring (see section~2 for a precise
definition).  We write $\GKdim$ for the Gelfand-Kirillov dimension.
The second major result of this paper is the following:
\begin{theorem} \label{purity}
Let $(\uu,\vv)$ be an integer vector which induces an increasing
filtration on $R$ such that $\gr_{(\uu,\vv)}(R)$ is a commutative
polynomial ring.  If $M$ is a finitely generated left $R$-module such
that, for every nonzero submodule $M'$ of $M$, one has $\GKdim M' \geq
p$, then every irreducible component of $\Ch_{(\uu,\vv)}(M)$ has
dimension at least $p$.
\end{theorem}
\noindent In particular, this theorem says that if $M$ is GKdim-pure,
that is $\GKdim M' = \GKdim M$ for every nonzero submodule $M'$ of
$M$, then the characteristic variety is equidimensional.  

Theorem~\ref{purity} generalizes known equidimensionality results to a
larger class of filtrations.  In particular, when $R$ is the
enveloping algebra of finite dimensional Lie algebra, it extends
Gabber's equidimensionality theorem~\cite[Th\'{e}or\`{e}me~1]{G-L}
beyond the standard filtration.  For certain skew polynomial rings, it
also extends the equidimensionality theorem in Huishi and van
Oystaeyen \cite[Corollary~III 4.3.6]{HvO} to non-Zariskian
filtrations.  Specifically, when the vector $(\uu,\vv)$ has a negative
entry, the induced filtration is not Zariskian.  Our proof of
Theorem~\ref{purity} involves studying the growth of filtered modules
and Gr\"{o}bner basis theory and differs significantly from the
homological methods used by Bj\"{o}rk~\cite{bj1}, Gabber~\cite{G-L}
and Huishi and van Oystaeyen~\cite{HvO}.

We now describe the contents of this paper.  In the next section, we
list the global notation and our general conventions.  In particular,
we recall the definitions of an almost centralizing extension, an
increasing filtration, an associated Rees ring and a characteristic
variety.  We also introduce the polynomial region, the space of
vectors inducing filtrations for which the associated graded algebra
is the commutative polynomial ring.  In Section~3, we use properties
of Rees modules to connect irreducible components of the
characteristic variety to the Gelfand-Kirillov dimension of
submodules.  To guarantee that the Gelfand-Kirillov dimension is well
behaved, we restrict our attention to finite dimensional filtrations
throughout this section.  Section~4 develops the Gr\"{o}bner basis
theory for our skew polynomial ring.  These tools are applied in
Section~5 to construct a combinatorial object, called the Gr\"{o}bner
fan.  This generalizes the Gr\"{o}bner fan of Mora and
Robbiano~\cite{MR} in the case of commutative polynomial rings and
Assi, Castro-Jim\'{e}nez and Granger~\cite{ACG} in the case of the
Weyl algebra.  In the last section, we use the use the Gr\"{o}bner fan
to extend our results for finite dimensional filtrations to all
filtrations in the polynomial region and prove our main theorems.

I would like to thank my advisor David Eisenbud for his support and
encouragement, Bernd Sturmfels for introducing me to algebraic
analysis, Nobuki Takayama for his help with his computer package {\tt
sm1/kan}, and Harrison Tsai for many useful discussions.

\section{Rings, Filtrations and Modules} 

In this section, we recall the notion of an almost centralizing
extension and list some of its basic properties.  Throughout this
paper $k$ denotes a field.

\subsection*{Almost centralizing extensions}
Let $B$ be the commutative polynomial ring $k[x_{1},\ldots,x_{m}]$.
We concentrate on a $k$-algebra $R$ which is generated by $x_{1},
\ldots, x_{m}, y_{1}, \ldots, y_{n}$ subject only to the relations:
\begin{align}
y_{i}x_{j} - x_{j}y_{i} & = Q^{1}_{i,j}(x) \, ; \tag{$R1$} \\
y_{i}y_{j}-y_{j}y_{i} & = Q^{2}_{i,j}(x,y) = Q^{2,0}_{i,j}(x) +
\sum_{\ell=1}^{n} Q^{2,\ell}_{i,j}(x) y_{\ell} \tag{$R2$}
\end{align}
where $Q^{1}_{i,j}(x)$, $Q^{2,\ell}_{i,j}(x) \in B$ for all $1 \leq
\ell \leq n$.  The skew polynomial ring $R$ is called an almost
centralizing extension of $B$.  The Poincar\'{e}-Birkhoff-Witt theorem
generalizes to $R$ and, hence, the set of standard monomial $\{
x^{\aa}y^{\bb} = x_{1}^{a_{1}} \cdots x_{m}^{a_{m}}y_{1}^{b_{1}}
\cdots y_{n}^{b_{1}} : (\aa,\bb) \in \NN^{m}\times\NN^{n} \}$ forms a
$k$-basis.  In particular, each element $f \in R$ has a unique {\sf
standard expression} of the form $\sum \kappa_{\aa,\bb}
x^{\aa}y^{\bb}$.  For more information about almost centralizing
extension see subsections~8.6.6 and 8.6.7 in McConnell and
Robson~\cite{McR}.

\begin{example}
If $\fg$ is a finite dimensional Lie algebra over $k$ then any crossed
product $B * U(\fg)$ is an almost centralizing extension of $B$.
Notably, the polynomial ring in $m+n$ central indeterminates over $k$,
the Weyl algebra $A_{n}(k)$ and the universal enveloping algebra
$U(\fg)$ all have this form.
\end{example}

\subsection*{The polynomial region}
A vector $(\uu,\vv) \in \RR^{m} \times \RR^{n}$ induces an increasing
filtration of $R$ as follows: for each $i \in \ZZ$, consider the
vector space
$
F_{i}R := k \cdot \left\{ x^{\aa}y^{\bb} : \text{$\lceil \uu \rceil
\cdot \aa + \lceil \vv \rceil \cdot \bb \leq i$, $\aa \in \NN^{m}$,
$\bb \in \NN^{n}$} \right\} \, ,
$
where $\lceil\uu\rceil = (\lceil u_{1} \rceil, \ldots, \lceil u_{m}
\rceil)$ and $\lceil u_{i}\rceil$ is the smallest integer greater than
$u_{i}$.  This clearly gives an increasing sequence of subspaces
satisfying the conditions $1 \in F_{0}R$ and $\bigcup_{i \in \ZZ}
F_{i}R = R$.  When $F_{i}R \cdot F_{j}R \subseteq F_{i+j}R$, the {\sf
associated graded ring} is $\gr_{(\uu,\vv)}(R) = \bigoplus_{i \in \ZZ}
F_{i}R / F_{i-1}R$.  An element $f$ belonging to the vector space
$F_{i}R - F_{i-1}R$ is said to have {\sf degree} $i$ and we write
$\deg_{(\uu,\vv)} f = i$.

The {\sf polynomial region} associated to $R$, denoted $\PR(R)$, is
set of all real vectors $(\uu,\vv)$ such that $\gr_{(\uu,\vv)}(R)$ is
the commutative polynomial ring generated by the initial forms
$\bar{x}_{i}$ and $\bar{y}_{i}$ of $x_{i}$ and $y_{i}$.  We denote
this commutative polynomial ring by $S$.  Since $S$ is noetherian
ring, $R$ is both a left and a right noetherian ring.  We will focus
exclusively on vectors $(\uu,\vv)$ belonging to the polynomial region.
The next proposition provides a more explicit interpretation of
$\PR(R)$.
\begin{proposition} 
The polynomial region $\PR(R)$ is the open convex polyhedral cone in
$\RR^{m}\times\RR^{n}$ given by the intersection of the
following open half-spaces:
\[ \begin{array}{rcll}
\deg_{(\uu,\vv)} x_{i}y_{j} & > & \deg_{(\uu,\vv)} Q^{1}_{i,j} &
\text{ for all $1 \leq j \leq n$ and $1 \leq i \leq m$}, \\
\deg_{(\uu,\vv)} y_{i}y_{j} & > & \deg_{(\uu,\vv)} Q^{2}_{i,j} &
\text{ for all $1 \leq i,j \leq n$.}
\end{array} \]
\end{proposition}

\begin{proof}
Without loss of generality, we may replace $(\uu,\vv)$ with $(\lceil
\uu \rceil,\lceil \vv \rceil)$.  Let $f_{1},\ldots,f_{j}$ be elements
of the $k$-vector space generated by the elements $x_{1}, \ldots,
x_{m}, y_{1}, \ldots, y_{n}$.  Now, if $\sigma$ is a permutation of
$\{1, \ldots, j \}$, we claim that
$
f_{1}f_{2} \cdots f_{j} \in f_{\sigma(1)}f_{\sigma(2)} \cdots
f_{\sigma(j)} + F_{\ell-1}R \, , 
$
where $\ell$ is the degree of the left hand side.  Since the elements
$x_{1}, \ldots, x_{m}$ commute and $Q_{i,j}^{1}(x)$ and
$Q_{i,j}^{2}(x,y)$ are at most linear in the variables
$y_{1},\ldots,y_{n}$, it suffices to prove this when $\sigma$ is a
transposition.  By linearity, the assertion is equivalent to the
conditions:
\begin{eqnarray*}
y_{i}x_{j} - x_{j}y_{i} & = & Q_{i,j}^{1}(x) \in F_{u_{i}+u_{j}-1}R \\
y_{i}y_{j} - y_{j}y_{i} & = & Q_{i,j}^{2}(x,y) \in F_{v_{i}+v_{j}-1}R
\, .
\end{eqnarray*}
We conclude that $F_{i}R \cdot F_{j}R \subseteq F_{i+j}R$ and that
$\gr_{(\uu,\vv)}(R)$ is a commutative $k$-algebra generated by
$\bar{x}_{1},\ldots,\bar{x}_{m}, \bar{y}_{1},\ldots, \bar{y}_{n}$.

To see that $\gr_{(\uu,\vv)}(R) = S$, it suffices to see that
there are no $k$-linear relations between the monomials in
$\gr_{(\uu,\vv)}(R)$.  A relation among the monomials in
$\gr_{(\uu,\vv)}(R)$, it would yield a relation among the standard
monomials in $R$.  However, the standard monomials form a $k$-basis
for $R$ which completes the proof.
\end{proof}

\begin{remark}
If $p = \max\{ \deg_{(\11,\00)} Q^{\ell}_{i,j} : \text{ for all $i$,
$j$, $\ell$} \}+1$, then the positive vector $(\11, p \11)$ belongs to
the polynomial region $\PR(R)$.
\end{remark}

\begin{example}
The polynomial region for $S$ is the entire space
$\RR^{m}\times\RR^{n}$ and $\PR(A_{n}(k)) = \{ (\uu,\vv) \in
\RR^{n}\times\RR^{n} : \text{$u_{i}+v_{i} > 0$ for all $1 \leq i \leq
n$}\}$. If $\mathfrak{sl}_{2}(k)$ has the standard basis $y_{1},
y_{2}, y_{3}$ such that $y_{2}y_{3}-y_{3}y_{2} = 2y_{3}$,
$y_{2}y_{1}-y_{1}y_{2} = -2y_{1}$ and $y_{1}y_{3}-y_{3}y_{1} = y_{2}$,
then $\PR(U(\mathfrak{sl}_{2}(k)))$ is the open cone in $\RR^{3}$ given by
the inequalities $v_{1}+v_{3} > v_{2}$ and $v_{2} > 0$.
\end{example}

\subsection*{Rees rings}
The {\sf associated Rees ring} of $R$ with respect to $(\uu,\vv)$ is
the graded $k$-module $\tR = \bigoplus_{i \in \ZZ} F_{i}R$.  The
$k$-algebra structure on $R$ makes $\tR$ into a graded $k$-algebra.
For $f \in F_{i}R$, we write the homogeneous element represented by
$f$ in $\tR_{i}$ as $(\tilde{f})_{i}$.  Observe that $(\tilde{1})_{1}
\in \tR_{1}$ is a central nonzero-divisor.  We denote this canonical
element by $x_{0}$.  More concretely, the Rees ring $\tR$ is generated
by $x_{0}, \ldots, x_{m}, y_{1}, \ldots, y_{n}$ subject to the
relations:
\begin{align}
y_{i}x_{j} - x_{j}y_{i} & = x_{0}^{\lceil u_{j} \rceil
+ \lceil v_{i} \rceil -q^{1}_{i,j}} Q^{1}_{i,j}(x); \tag{$\tR1$} \\
y_{i}y_{j}-y_{j}y_{i} & = x_{0}^{\lceil v_{i} \rceil +
\lceil v_{j} \rceil -q^{2}_{i,j}} Q^{2}_{i,j}(x,y), \tag{$\tR2$}
\end{align}
where $Q^{1}_{0,j}(x) = 0$ for all $1 \leq j \leq n$ and
$q^{\ell}_{i,j} = \deg_{(\uu,\vv)} Q^{\ell}_{i,j}$ for $1 \leq \ell
\leq 2$.  We stress that $\tR$ is an almost centralizing extension of
$B[x_{0}]$ and the relations ($\tR$1) and ($\tR$2) are homogeneous
with respect to $(\uu,\vv)$.  The condition that $(\uu,\vv)$ belongs
to the polynomial region $\PR(R)$ insures that $x_{0}$ has a
nonnegative exponent in relations ($\tR$1) and ($\tR$2).  The
homogenization map from $R$ to $\tR$ is defined as follows: for $f =
\sum_{\aa,\bb} \kappa_{\aa,\bb} x^{\aa}y^{\bb}$ in $R$, we set
$\tilde{f} = \sum_{\aa,\bb} \kappa_{\aa,\bb} x_{0}^{i -
\uu\cdot\aa-\vv\cdot\bb}x^{\aa}y^{\bb}$ where $\deg_{(\uu,\vv)}(f)=i$.
In the other direction, the substitution $x_{0}=1$ gives a $k$-algebra
homomorphism $\tR$ to $R$.

\subsection*{Filtered modules and characteristic varieties}
All modules and ideals considered in this paper will be finitely
generated left modules and left ideals respectively.  We write $M$ for
a finitely generated $R$-module.  By a filtered $R$-module, we mean
that there is an increasing sequence of $k$-vector spaces $F_{i}M$ for
$i \in \ZZ$ satisfying the conditions: $F_{i}R \cdot F_{j}M \subseteq
F_{i+j}M$, and $\bigcup_{i \in \ZZ} F_{i}M = M$.  We define the {\sf
associated graded module} to be $\gr(M) = \bigoplus_{i \in \ZZ} F_{i}M
/ F_{i-1} M$.  It follows that $\gr(M)$ is a graded
$\gr_{(\uu,\vv)}(R)$-module.  

A {\sf good filtration} is a filtration of an $R$-module $M$ for which
there exists elements $z_{1},\ldots,z_{p}$ in $M$ and integers
$w_{1},\ldots,w_{p}$ such that $F_{i}M = \sum_{j=1}^{p} F_{i-w_{j}}R
\cdot z_{j}$.  Every finitely generated $R$-module $M$ has a good
filtration and, conversely, any module with a good filtration is
necessarily finitely generated over $R$.  For a good filtration of
$M$, we define the {\sf characteristic ideal} $I(M)$ to be the prime
radical of $\Ann_{S}(\gr(M))$.  Since any two good filtrations
are equivalent, the characteristic ideal $I(M)$ is independent of the
choice of good filtration; however $I(M)$ does depend on $(\uu,\vv)$.
The {\sf characteristic variety} of $M$ is the reduced scheme
$\Ch_{(\uu,\vv)}(M) = \Spec( \gr_{(\uu,\vv)}(R) / I(M))$.

For a filtered $R$-module $M$, we define the {\sf associated Rees
module} to be the graded $k$-module $\tM = \bigoplus_{i \in \ZZ}
F_{i}M$.  The $R$-module structure on $M$ makes $\tM$ into a graded
$\tR$-algebra.  More details on Rees modules can be found in
section~I.4 of Huishi and van Oystaeyen~\cite{HvO}.

\section{Finite Dimensional Filtrations}

Under the assumption that $R$ has a finite dimensional filtration, we
are able to relate the dimension of the irreducible components of
$\Ch_{(\uu,\vv)}(M)$ to the Gelfand-Kirillov dimension of submodules
of $M$.  We accomplish this by using the Rees module $\tM$ to link
submodules of $M$ and graded submodules of $\gr(M)$.  We begin with a
brief discussion of Gelfand-Kirillov dimension.

We define the Gelfand-Kirillov dimension only for $R$-modules with a
given finite dimensional filtration.  Recall that a function $\phi
\colon \NN \rightarrow \RR_{\geq 0}$ has {\sf polynomial growth}, if
for some $d \in \RR$, $\phi(i) \leq i^{d}$ for all $i$ sufficiently
large.  In this situation, we consider the number $\gamma(\phi) :=
\inf\left\{ d : \text{$f(i) \leq i^{d}$ for $i \gg 0$} \right\}$.  For
a filtered $R$-module $M$, the {\sf Gelfand-Kirillov
dimension} is
$
\GKdim M := \gamma\left( \dim_{k} F_{i}M \right) \, .  
$
It follows, from subsection~8.6.18 in McConnell and Robson~\cite{McR},
that $\GKdim M$ is independent of choice of generators $z_{j}$ and
integers $w_{j}$, although it depends on the filtration of $R$ --- see
Proposition~\ref{GKdim-indep} for a discussion of this dependence.

Now, if $(\uu,\vv) \in \PR(R)$ is not positive, then $\dim_{k}F_{i}M$
is infinite for all $i$.  Thus, for $R$ to have a finite dimensional
filtration, it is necessary and sufficient that the vector $(\uu,\vv)$
be positive.  With this additional hypothesis, we can give a useful
description of the function $i \mapsto \dim_{k} F_{i}M$.  Recall that
a function $\theta \colon \ZZ \rightarrow \CC$ is called a {\sf
quasi-polynomial} if there exists a positive integer $p$ and
polynomials $Q_{j}$ for $0 \leq j \leq p-1$ such that, for all $i \in
\ZZ$, we have $\theta(i) = Q_{j}(i)$ where $i = rp+ j$ with $0 \leq j
\leq p-1$.  The degree of a quasi-polynomial is the maximum of degree
of the polynomials $Q_{j}$.  The next proposition is a small extension
of the almost commutative results in section~8.4 of McConnell and
Robson~\cite{McR}.

\begin{proposition} \label{hilbert}
Assume the vector $(\uu,\vv) \in \PR(R)$ is positive.  If $M$ is a
nonzero $R$-module with a good finite dimensional filtration such that
$\gr(M)$ has Krull dimension $d$, then one has the following:
\begin{enumerate}
\item[(1)] There exists positive integers $c_{0}, \ldots,
c_{d}$ and $Q(t) \in \ZZ[t,t^{-1}]$ such that
\[
\sum_{i \geq 0} \big( \dim_{k} F_{i}M \big) \cdot t^{i}
= \frac{ Q(t)}{ \prod_{j=0}^{d} (1-t^{c_{j}})} \, ,
\text{  with $Q(1) > 0$.}
\]
\item[(2)] The function $i \mapsto \dim_{k} F_{i}M$ is a
quasi-polynomial of degree $d$.
\end{enumerate}
\end{proposition}

\begin{proof}
Since $(\uu,\vv) \in \PR(R)$ is positive, $S$ is positively graded
commutative $k$-algebra and Proposition~4.4.1 in Bruns and
Herzog~\cite{BH} implies there are positive integers $c_{1}, \ldots,
c_{d}$ and $Q(t) \in \ZZ[t,t^{-1}]$ such that $\sum_{i \geq 0}
\big( \dim_{k} \gr(M)_{i} \big) \cdot t^{i} = Q(t) /
\prod_{j=1}^{d} (1-t^{c_{j}})$ and $Q(1) > 0$.  Hence, we have
\begin{eqnarray*}
\sum_{i \geq 0} \big( \dim_{k} F_{i}M \big) \cdot t^{i}  
= 
\sum_{i \geq 0} \left( \sum_{j=0}^{i} \dim_{k} \gr(M)_{j} \right)
\cdot t^{i} \\ 
=  
\left( \frac{ Q(t)}{
\prod_{j=1}^{d} (1-t^{c_{j}})} \right) \left( \frac{1}{1-t}
\right),
\end{eqnarray*}
which proves part~(1).  Part~(2) follows immediately from part~(1) by
applying Proposition~4.4.1 in Stanley~\cite{Sta}.
\end{proof}

The second part of this proposition clearly implies the following:

\begin{corollary} \label{commutative}
If the vector $(\uu,\vv) \in \PR(R)$ is positive and $N$ is a finitely
generated graded $S$-module, then the Gelfand-Kirillov dimension and
Krull dimension of $N$ are equal.\qed
\end{corollary}

We now turn our attention to Rees modules and provide homomorphisms
linking the submodules of $M$, $\tM$ and $\gr(M)$.  We always assume
that the filtration on a submodule is the unique filtration induced by
the module containing it.  Because $R$ is left noetherian, good
filtrations induce good filtrations on submodules.

\begin{proposition} \label{pi}
Let $x_{0}$ be the canonical homogeneous element of degree $1$ in Rees
ring $\tR$.  If $M$ is a filtered $R$-module and $\tM$ is the
associated Rees module, one has the following:
\begin{enumerate}
\item[(1)] There exists a surjective homomorphism $\pi_{1} \colon \tM
\to M$ such that $\Ker \pi_{1} = (1-x_{0}) \cdot M$.  Moreover, for all
submodules $M' \subseteq M$, one has $\pi_{1}(\widetilde{M'}) = M'$.
\item[(2)] There exists a surjective graded homomorphism $\pi_{0}
\colon \tM \to \gr(M)$ such that $\Ker \pi_{0} = x_{0} \cdot M$.
Furthermore, $\pi_{0}$ maps graded submodules of $\tM$ to graded
submodules of $\gr(M)$ and every graded submodule of $\gr(M)$ arises
in this manner.
\end{enumerate}
\end{proposition}

\begin{proof}
(1) Every element $\tilde{z} \in \tM$ can be written uniquely as a
finite sum of homogeneous components; $\tilde{z} = \sum_{j=0}^{p}
(\tilde{z})_{\ell_{j}}$ where $\ell_{0} < \cdots < \ell_{p}$.  Let
$\pi_{1} \colon \tM \to M$ be defined by $\pi_{1}(\tilde{z}) =
\sum_{j=0}^{p} (z)_{\ell_{j}}$ where $(z)_{\ell_{j}} \in
F_{\ell_{j}}M$.  The definition of the $\tR$-module structure on $\tM$
insures that $\pi_{1}$ is a $k$-module homomorphism and the image is
an $R$-module.  It is clearly surjective.  Now, if $\sum_{j=0}^{p}
(z)_{\ell_{j}} = 0$ then $\sum_{j=0}^{p}
(\tilde{z})_{\ell_{j}}x_{0}^{\ell_{p}-\ell_{j}} = 0$.  Hence, the
element
\[
\tilde{z} = \sum_{j=0}^{p} (\tilde{z})_{\ell_{j}} -
\sum_{j=0}^{p} (\tilde{z})_{\ell_{j}}x_{0}^{\ell_{p}-\ell_{j}}
= \sum_{j=0}^{p} \big( (\tilde{z})_{\ell_{j}} -
(\tilde{z})_{\ell_{j}}x_{0}^{\ell_{p}-\ell_{j}} \big)
\]
belongs to $(1-x_{0})\cdot\tM$ and we have $\Ker \pi_{1} \subseteq
(1-x_{0})\cdot\tM$.  It is obvious from the definition of $\pi_{1}$ that
we have $(1-x_{0})\cdot\tM \subseteq \Ker \pi_{1}$ and
$\pi_{1}(\widetilde{M'}) = M'$.

(2) For all $i \in \ZZ$, we have isomorphisms 
\[
\tM_{i} / (x_{0}\cdot\tM_{i-1}) \cong F_{i}M / F_{i-1}M = \gr(M)_{i}
\, .
\] 
Combining this maps gives the required isomorphism $\tM /
(x_{0}\cdot\tM) \cong \gr(M)$.  Moreover, we have
\[
\pi_{0}(\tilde{f}\tilde{z}) = fz+(x_{0}\cdot M) = (f+x_{0}\cdot R)(z +
x_{0}\cdot M) = \pi_{0}(\tilde{f})\pi_{0}(\tilde{z})
\]
and, thus, $\pi_{0}$ takes $\tR$-modules to $\gr(R)$-modules.  Finally,
for a graded submodule $N$ of $\gr(M)$, consider the $\tR$-submodule
$L$ of $\tM$ generated by the set $\pi_{0}^{-1}(N)$.  To demonstrate
that $\pi_{0}(L) = N$, it suffices to show $\pi_{0}(L) \subseteq N$.
Every element of $L$ can be written in the form $\sum_{j=0}^{p}
\tilde{f}_{j} \tilde{z}_{j}$ for some $\tilde{f}_{j} \in \tR$ and
$\tilde{z}_{j} \in \pi_{0}^{-1}(N)$.  Applying $\pi_{0}$, we obtain
$$ \pi_{0}\left(\sum_{j=0}^{p} \tilde{f}_{j} \tilde{z}_{j} \right)
= \sum_{j=0}^{p} \pi_{0}(\tilde{f}_{j}) \pi_{0}(\tilde{z}_{j})
= \sum_{j=0}^{p} f_{j} z_{j} $$
where $f_{j} \in \gr(A)$ and $z_{j} \in N$.  Therefore,
we have $\pi_{0}(L) \subseteq N$ which completes the proof.
\end{proof}

We next record a useful lemma; see Proposition~2.16 in
Bj\"{o}rk~\cite{bj2}.

\begin{lemma}[Bj\"{o}rk] \label{torsion}
Let $M$ be a filtered $R$-module.  If $L$ be a graded submodule of
$\tM$, then the graded module $\widetilde{\pi_{1}(L)}$ contains $L$
and the quotient $\widetilde{\pi_{1}(L)} / L$ is an $x_{0}$-torsion
module. \qed
\end{lemma}

\begin{proposition} \label{pi1dim}
If $(\uu,\vv) \in \PR(R)$ is positive, $M$ is an $R$-module with a
good finite dimensional filtration and $L$ is a graded submodule of
$\tM$, then one has the following:
\begin{enumerate}
\item[(1)] $\GKdim L = \GKdim \widetilde{\pi_{1}(L)}$;
\item[(2)] $1 + \GKdim M = \GKdim \tM$;
\item[(3)] $1 + \GKdim \pi_{1}(L) = \GKdim L$.
\end{enumerate}
\end{proposition}

\begin{proof}
(1) Let $\phi(i) = \dim_{k} L_{i}$, $L' := \widetilde{\pi_{1}(L)}$ and
$\psi(i) = \dim_{k} L_{i}'$.  By Lemma~\ref{torsion}, $L_{i}$ is a
subvector space of $L_{i}'$ which implies $\phi(i) \leq \psi(i)$ and
$\GKdim L \leq \GKdim L'$.  On the other hand, Lemma~\ref{torsion}
also states that the quotient $L' / L$ is an $x_{0}$-torsion module.
Since $L'$ is a finitely generated module, there exists an integer
$\ell$ such that $x_{0}^{\ell} \cdot L_{i}' \subseteq L_{i+\ell}$.
Thus, we have $\psi(i) \leq \phi(i+\ell)$ which implies $\GKdim L'
\leq \GKdim L$.  Combining the two inequality yields the first part.

(2) The definition of Gelfand-Kirillov dimension implies
\begin{eqnarray*}
\GKdim M & = & \gamma\left( \dim_{k} F_{i}M \right) \text{ and } \\
\GKdim \tM & = & \gamma\left( \sum_{j=0}^{i} \dim_{k} F_{j}M \right)
\, .
\end{eqnarray*}
However, a monotonically increasing function $\phi \colon \NN \to
\RR_{\geq 0}$ and the function $\psi(i) = \sum_{j=0}^{i} \phi(j)$ are
related by the equation $\gamma(\psi) = \gamma(\phi)+1$ and this
proves the second assertion.

(3) Applying part~(2) gives $1 + \GKdim \pi_{1}(L) = \GKdim
\widetilde{\pi_{1}(L)}$ and combining this part~(1) yields the third
assertion.
\end{proof}

We have the analogous result for submodules of $\tM$ and $\gr(M)$.

\begin{proposition} \label{pi0dim}
Let $(\uu,\vv) \in \PR(R)$ be a positive vector and let $M$ be an
$R$-module with a good filtration.  If $L$ is a graded $\tR$-submodule
of $\tM$ then $1 + \GKdim \pi_{0}(L) = \GKdim L$.
\end{proposition}

\begin{proof}
We begin by observing the following:
\[
\dim_{k} F_{i}\big(\pi_{0}(L)\big)
= \sum_{j=0}^{i} \dim_{k} \frac{ L_{j}}{x_{0} \cdot L_{j-1}}
= \dim_{k} L_{i} \, .
\]
Next, notice that $\dim_{k} F_{i}L = \sum_{j=0}^{i} \dim_{k} L_{j}$.
Since $x_{0}$ is a nonzero-divisor of degree $1$ on $L$, we have
$\dim_{k} L_{i} \leq \dim_{k} L_{i+1}$.  Thus, applying the formula
$\gamma(\psi) = \gamma(\phi)+1$ for a monotonically increasing
function $\phi \colon \NN \to \RR_{\geq 0}$ and the function $\psi(i)
= \sum_{j=0}^{i} \phi(j)$, completes the proof.
\end{proof}

We now link the dimension of the irreducible components of the
characteristic variety $\Ch_{(\uu,\vv)}(M)$ to the Gelfand-Kirillov
dimension of submodules of $M$.

\begin{theorem} \label{GKdim-Kdim}
Let $p$ be a nonnegative integer, let $(\uu,\vv) \in \PR(R)$ be a
positive vector and let $M$ be a filtered $R$-module.  If
$\Ch_{(\uu,\vv)}(M)$ has an irreducible component of dimension $p$
then there exists a submodule $M'$ of $M$ such that $\GKdim M' = p$.
\end{theorem}

\begin{proof}
By definition, irreducible components of $\Ch_{(\uu,\vv)}(M)$
correspond to minimal primes in the support of $\gr(M)$.  Hence, if
there exists an irreducible component of $\Ch_{(\uu,\vv)}(M)$ of
dimension $p$, we have a minimal prime $\fp$ in the support of
$\gr(M)$ of dimension $p$.  Each minimal prime $\fp$ in the support of
$\gr(M)$ corresponds to a graded submodule of $\gr(M)$ of the form
$\big( S / \fp \big)(j)$ for some $j \in \ZZ$ and
Corollary~\ref{commutative} implies the Krull dimension of $\big( S /
\fp \big)(j)$ is equal to its Gelfand-Kirillov dimension.  Thus, we
have a graded submodule of $\gr(M)$ with Gelfand-Kirillov dimension
$p$.  We complete this proof by showing that the following three
conditions are equivalent:
\begin{enumerate}
\item[(a)] there exist a submodule $M'$ of $M$ with $\GKdim M' =
p$;
\item[(b)] there exist a graded submodule $L$ of $\tM$ with
$\GKdim L = p+1$;
\item[(c)] there exist a graded submodule $N$ of $\gr(M)$ with
$\GKdim N = p$;
\end{enumerate}
Indeed, we have:
\begin{enumerate}
\item[(a)] $\Rightarrow$ (b): By Proposition~\ref{pi1dim}.2, the
graded submodule $\widetilde{M'}$ of $\tM$ has Gelfand-Kirillov
dimension $p+1$.
\item[(b)] $\Rightarrow$ (a): Follows immediately from
Proposition~\ref{pi1dim}.3.
\item[(b)] $\Rightarrow$ (c): Follows immediately from
Proposition~\ref{pi0dim}.
\item[(c)] $\Rightarrow$ (b): By Proposition~\ref{pi}, there exists a
graded submodule $L$ of $\tM$ such that $\pi_{0}(L)=N$.
\end{enumerate}
\end{proof}

\begin{remark} \label{upper_bound}
Since $\GKdim M = \GKdim \gr(M)$ for any $R$-module $M$, we see that
the Gelfand-Kirillov dimension of $M$ is an upper bound on the Krull
dimension of each irreducible component of $\Ch_{(\uu,\vv)}(M)$.
\end{remark}

We point out that Theorem~\ref{GKdim-Kdim} can also be proven using
homological methods; see Chapter~2 in Bj\"{o}rk \cite{bj1} or
Chapter~III in Huishi and van Oystaeyen \cite{HvO} for the techniques.
However, the approach presented here is more elementary.

\begin{remark}
We have not used the fact that $R$ is an almost centralizing extension
of $B$, so Proposition~\ref{pi}, Proposition~\ref{pi1dim},
Proposition~\ref{pi0dim} and Theorem~\ref{GKdim-Kdim} hold a module
$M$ with a good finite dimensional filtration over a filtered
$k$-algebra.
\end{remark}

\section{Gr\"{o}bner Basics}

This section is devoted to the Gr\"{o}bner basics of the skew
polynomial ring $R$.  Sturmfels' expression ``Gr\"{o}bner basics''
describes a collection of ideas centering around initial ideals with
respect to a vector $(\uu,\vv)$ and term orders.  We develop the
algebraic aspects of this theory, generalizing results in commutative
polynomial ring and Weyl algebra.  

\subsection*{Gr\"{o}bner basis with respect to $(\uu,\vv)$}
For $f = \sum \kappa_{\aa,\bb} x^{\aa}y^{\bb}$ in $R$, the {\sf
initial form} (also called the {\sf principal symbol}) of $f$ with
respect to $(\uu,\vv)$ is the element $\IN_{(\uu,\vv)}(f) =
\sum_{\uu\cdot\aa + \vv\cdot\bb = \ell} \kappa_{\aa,\bb}
\bar{x}^{\aa}\bar{y}^{\bb}$ in $S$, where $\ell = \max\{
\uu\cdot\aa+\vv\cdot\bb : \text{$ \kappa_{\aa,\bb} \neq 0$} \}$.

\begin{proposition} \label{in=gr}
If $I$ is an $R$-ideal and $(\uu,\vv) \in\PR(R)$ then
\[
\IN_{(\uu,\vv)}(I) = k \cdot \left\{ \IN_{(\uu,\vv)}(f) : \text{$f \in
I$} \right\}
\] 
is an $S$-ideal.  Moreover, if $(\uu,\vv)$ is an integer vector
belonging to $\PR(R)$ then $\gr(I)$ is isomorphic to
$\IN_{(\uu,\vv)}(I)$.
\end{proposition}

\begin{proof}
By definition, $\IN_{(\uu,\vv)}(I)$ is closed under addition, so it is
enough to show that it is closed under left multiplication by the
elements $\bar{x}_{1}, \ldots, \bar{x}_{m}, \bar{y}_{1}, \ldots,
\bar{y}_{n}$.  Consider the element $f = \sum \kappa_{\aa,\bb} x^{\aa}
y^{\bb}$ in $I$.  One clearly has
\[
\bar{x}_{i} \cdot \IN_{(\uu,\vv)}(f) = \IN_{(\uu,\vv)}\left( \sum
\kappa_{\aa,\bb} x^{\aa+\ee_{i}} y^{\bb} \right) =
\IN_{(\uu,\vv)}(x_{i} \cdot f) \, ,
\]
where $\ee_{i}$ is the $i$-th standard basis vector.  Similarly,
because $(\uu,\vv)$ belongs to $\PR(R)$, we obtain
\[
\bar{y}_{i} \cdot \IN_{(\uu,\vv)}(f) = \IN_{(\uu,\vv)}\left( \sum
\kappa_{\aa,\bb} x^{\aa} y^{\bb+\ee_{i}} + Q(x,y) \right) =
\IN_{(\uu,\vv)}(y_{i} \cdot f) \, ,
\]
where $Q(x,y) \in R$ has an initial form which is strictly smaller
than $\IN_{(\uu,\vv)}\left(\sum \kappa_{\aa,\bb} x^{\aa}
y^{\bb+\ee_{i}} \right)$.

When $(\uu,\vv) \in \PR(R) \cap \ZZ^{m} \times \ZZ^{n}$, we note that
\[
\{ \IN_{(\uu,\vv)}(f) : \text{ $f \in I$ and $\deg_{(\uu,\vv)}(f) =
i$}\}
\] 
is a complete set of representatives for the cosets of $\gr(I)_{i}$.
Hence, there exists a bijective set map between $\IN_{(\uu,\vv)}(I)$
and $\gr(I)$ and one easily verifies that the $S$-module structure of
$\IN_{(\uu,\vv)}(I)$ and $\gr(I)$ agree under this correspondence.
\end{proof}

The $S$-ideal $\IN_{(\uu,\vv)}(I)$ is called the {\sf initial ideal}
of the $R$-ideal $I$ with respect to the vector $(\uu,\vv)$.  A finite
subset $\cG$ of $R$ is a {\sf Gr\"{o}bner basis} of $I$ with respect
to $(\uu,\vv)$ if $I$ is generated by $\cG$ and $\IN_{(\uu,\vv)}(I)$
is generated by the initial forms $\IN_{(\uu,\vv)}(\cG) = \left\{
\IN_{(\uu,\vv)}(g) : g \in \cG \right\}$.

\subsection*{Gr\"{o}bner bases with respect to $\prec$}
A total ordering $\prec$ on the set of standard monomials $\left\{
x^{\aa}y^{\bb} : (\aa,\bb) \in \NN^{m}\times\NN^{n} \right\}$ in $R$
is called a {\sf multiplicative order} when the following three
conditions hold:
\begin{align}
& x^{\aa} \prec x_{i}y_{j}  \text{ for all monomials
$x^{\aa}$ appearing in $Q^{1}_{i,j}(x)$;} \tag{M1} \\
& x^{\aa}y_{\ell} \prec y_{i}y_{j} \text{ for all monomials
$x^{\aa}y_{\ell}$ appearing in $Q^{2}_{i,j}(x,y)$;} \tag{M2} \\
& x^{\aa}y^{\bb} \prec x^{\aa'}y^{\bb'} \Rightarrow
x^{\aa+\cc}y^{\bb+\dd} \prec x^{\aa'+\cc} y^{\bb'+\dd} \,\, \forall \,
(\cc,\dd) \in \NN^{m} \times \NN^{n} \,. \tag{M3}
\end{align}
A multiplicative order $\prec$ is called a {\sf term order} if $1 =
x^{0}y^{0}$ is the smallest element of $\prec$.  A multiplicative
order which is not a term order has infinite strictly decreasing
chains but a term order does not.  For information on frequently used
term orders such as lexicographic order or reverse lexicographic
order, see chapter~15 in Eisenbud~\cite{E}.

\begin{remark}
Conditions (M1) and (M2) correspond directly to the relations ($R1$)
and ($R2$) in the definition of $R$.  Without these assumptions the
order would not be compatible with multiplication; that is we would
not have $\IN_{\prec}(f \cdot f') = \IN_{\prec}(f)\IN_{\prec}(f')$.
\end{remark}

Fix a multiplicative order $\prec$.  The {\sf initial monomial} of $f
\in R$ is the monomial $\bar{x}^{\aa}\bar{y}^{\bb} \in S$ such that
$x^{\aa}y^{\bb}$ is the $\prec$-largest monomial appearing in the
standard expansion of $f$ in $R$.  For an $R$-ideal $I$, the {\sf
initial ideal} $\IN_{\prec}(I)$ is the monomial ideal in $S$ generated
by $\left\{ \IN_{\prec}(f) : \text{ $f \in I$} \right\}$.  A finite
subset $\cG$ of $R$ is a {\sf Gr\"{o}bner basis} of $I$ with respect
to $\prec$ if $I$ is generated by $\cG$ and $\IN_{\prec}(I)$ is
generated by $\IN_{\prec}(\cG) = \left\{ \IN_{\prec}(g) : g \in \cG
\right\}$.  The Gr\"{o}bner basis is called {\sf reduced} if, for any
two distinct elements $g$, $g' \in \cG$, the exponent vector of
$\IN_{\prec}(g)$ is componentwise larger than any exponent vector
appearing in the normally order expression of $g'$ in $R$.

\subsection*{Comparing $\prec$ and $(\uu,\vv)$}
We have two different notions of Gr\"{o}bner basis in $R$; one for
degree vectors $(\uu,\vv)$ and one for multiplicative monomial orders.
Proposition~\ref{prec=>(u,v)} below relates these two notions.
Before giving this relation, we collect some preliminary results.  Let
$(\uu,\vv) \in \PR(R)$ and let $\prec$ be any term order.  The
multiplicative monomial order $\prec_{(\uu,\vv)}$ is defined as
follows:
\[ \begin{array}{rl}
x^{\aa'}y^{\bb'} \prec x^{\aa}y^{\bb} & \text{ if and only if } \\ 
& \big(
\uu\cdot(\aa-\aa')+\vv\cdot(\bb-\bb') > 0 \big) \\ 
& \text{ or } \left( \begin{array}{l}
\uu\cdot(\aa-\aa')+\vv\cdot(\bb-\bb') = 0 \\ 
\text{ and } x^{\aa'}y^{\bb'} \prec
x^{\aa}y^{\bb}
\end{array} \right).
\end{array} \]
Note that $\prec_{(\uu,\vv)}$ is a term order if and only if
$(\uu,\vv)$ is a non-negative vector.

\begin{lemma} \label{combo-term-order}
If $I$ is an $R$-ideal and $(\uu,\vv) \in \PR(R)$, then 
\[ 
\IN_{\prec}\left(\IN_{(\uu,\vv)}(I)\right) =
\IN_{\prec_{(\uu,\vv)}}(I) \, .
\]
\end{lemma}

\begin{proof}
See Proposition~1.8 is Sturmfels~\cite{Stu}.
\end{proof}

\begin{proposition} \label{prec=>(u,v)}
Let $I$ be any $R$-ideal, $(\uu,\vv) \in \PR(R)$ and let $\prec$ be
any term order.  If $\cG$ a Gr\"{o}bner basis for $I$ with respect to
$\prec_{(\uu,\vv)}$, then one has
\begin{enumerate}
\item[(1)] the set $\cG$ is a Gr\"{o}bner basis for $I$ with respect
to $(\uu,\vv)$; 
\item[(2)] the set $\IN_{(\uu,\vv)}(\cG) =\left\{ \IN_{(\uu,\vv)}(g)
: g \in \cG \right\}$ is a Gr\"{o}bner basis for $\IN_{(\uu,\vv)}(I)$
with respect to $\prec$;  
\item[(3)] if $\cG$ is the reduced Gr\"{o}bner basis for $I$ with
respect to $\prec_{(\uu,\vv)}$, then $\IN_{(\uu,\vv)}(\cG)$ is also
reduced.
\end{enumerate}
\end{proposition}

\begin{proof}
Parts (1) and (2) are analogous to Theorem~1.1.6 in Saito, Sturmfels
and Takayama~\cite{SST}.  Part (3) follows from the fact that the
exponent vectors appearing in $\IN_{(\uu,\vv)}(g)$ for $g \in \cG$
form a subset of the exponent vectors appearing the standard
expression of $g$ in $R$.
\end{proof}

\subsection*{The number of initial ideals}
Although there are infinitely many different term orders when $m+n
\geq 2$, this does not lead to an infinite number of distinct initial
ideals.  

\begin{theorem} \label{finite}
An $R$-ideal $I$ has only finitely many distinct initial ideals
$\IN_{\prec}(I)$ where $\prec$ is a term order.
\end{theorem}

\begin{proof}
See Theorem~1.2 in Sturmfels~\cite{Stu}.
\end{proof}

\subsection*{Homogenization}  
A multiplicative order $\prec$ on $R$ lifts to a multiplicative order
$<$ on $\tR$ by the following convention:
\[ \begin{array}{rl}
x_{0}^{a_{0}'}x^{\aa'}y^{\bb'} <
x_{0}^{a_{0}}x^{\aa}y^{\bb} & \text{ if and only if } \\ 
& \big( a_{0}' - a_{0} > 0 \big) \\ 
& \text{ or } \left( \begin{array}{l}
a_{0}' - a_{0} = 0 \\ 
\text{ and } x^{\aa'}y^{\bb'} \prec x^{\aa}y^{\bb}
\end{array} \right).
\end{array} \]
Note that $\prec$ is a term order if and only if $<$ is a term order.

\begin{proposition} \label{nontermorder}
Let $\prec$ be a multiplicative order on $R$ and let $\tI$ be the
homogenization of an $R$-ideal $I$ with respect to $(\uu,\vv) \in
\PR(R)$.  If $\tG$ is a Gr\"{o}bner basis for $\tI$ with respect to
$<$ then its dehomogenization $\cG$ is a Gr\"{o}bner basis for $I$
with respect to $\prec$.
\end{proposition}

\begin{proof}
Since $\tI|_{x_{0}=1} = I$, the set $\cG$ generates $I$ if and only if
$\tG$ generates $\tI$.  Thus, it suffices to study the initial ideals.
If $h \in \IN_{\prec}(I)$ then we have $\tilde{h} \in \IN_{<}(\tI)$.
Since $\tG$ is a Gr\"{o}bner basis, it follows that $\tilde{h} =
\IN_{<}(\tilde{f} \tilde{g})$ for some $\tilde{f} \in \tR$ and
$\tilde{g} \in \tG$.  Dehomogenizing, we obtain $h =
\tilde{h}|_{x_{0}=1} = \IN_{\prec}( \tilde{f}|_{x_{0}=1} \cdot
\tilde{g}|_{x_{0}=1})$ which implies $h \in \IN_{\prec}(\cG)$ as
required.
\end{proof}

Fix an $R$-ideal $I$.  Two degree vectors $(\uu,\vv)$ and
$(\uu',\vv')$ in $\PR(R)$ are {\sf equivalent with respect to $I$} if
$\IN_{(\uu,\vv)}(I) = \IN_{(\uu',\vv')}(I)$.  We denote the
equivalence class of vectors $(\uu,\vv)$ with respect to $I$ by
$C_{I}[(\uu,\vv)]$.  We define the {\sf Gr\"{o}bner region} $\GR(I)$
to be the set of all $(\uu,\vv) \in \PR(R)$ such that
$\IN_{(\uu,\vv)}(I) = \IN_{(\uu',\vv')}(I)$ for some positive vector
$(\uu',\vv') \in \PR(R)$.  

\begin{proposition} \label{GR=PR}
Suppose that $R$ is a graded $k$-algebra with respect to a positive
vector $(\uu,\vv) \in \PR(R)$.  If $I$ is a homogeneous $R$-ideal,
then we have $\GR(I) = \PR(R)$.
\end{proposition}

Notice that, for any $R$-ideal $I$, the Rees ring $\tR$ associated to
a positive vector $(\uu,\vv) \in \PR(R)$ and $\tI$ satisfy the
hypothesis.

\begin{proof}
See Proposition~1.12 in Sturmfels~\cite{Stu}.
\end{proof}

A finite subset $\cU$ of an $R$-ideal $I$ is called a {\sf universal
Gr\"{o}bner basis} if $\cU$ is simultaneously a Gr\"{o}bner basis of
$I$ with respect to all vectors $(\uu,\vv) \in \PR(R)$.  This
definition is different than the one found Sturmfels~\cite{Stu};
Sturmfels' considers only vectors $(\uu,\vv)$ in the Gr\"{o}bner
region $\GR(I)$.  Proposition~\ref{GR=PR} shows that these two
different notions of a universal Gr\"{o}bner basis agree for
homogeneous ideals in a graded ring.

\begin{corollary}
Every $R$-ideal $I$ has a finite universal Gr\"{o}bner basis.
\end{corollary}

\begin{proof}
Consider homogenization $\tI$ of $I$ with respect to a positive vector
$(\uu,\vv) \in \PR(R)$.  By Theorem~\ref{finite} there exists only
finitely many distinct reduced Gr\"{o}bner basis for $\tI$ with
respect to term orders.  Let $\tG$ be their union.  Fix a term order
$\prec$ on $R$ and let $<_{(\uu,\vv)}$ denote the multiplicative order
on $\tR$ obtained from the multiplicative order $\prec_{(\uu,\vv)}$ on
$R$.  By construction, $\tG$ is a Gr\"{o}bner basis with respect to
$<_{(\uu,\vv)}$ when $(\uu,\vv) \in \PR(R)$ is positive.  Applying
Proposition~\ref{combo-term-order} and Proposition~\ref{GR=PR}, we see
that $\tG$ is, in fact, a Gr\"{o}bner basis with respect to
$<_{(\uu,\vv)}$ when $(\uu,\vv) \in \PR(R)$.  If $\cG$ is the
dehomogenization of $\tG$, then Proposition~\ref{nontermorder} implies
that $\cG$ is a Gr\"{o}bner basis of $I$ for $\prec_{(\uu,\vv)}$ where
$(\uu,\vv) \in \PR(R)$.  Finally, Proposition~\ref{prec=>(u,v)}.1
shows that $\cG$ is a universal Gr\"{o}bner basis for $I$.
\end{proof}

\begin{remark}
Since they are not required for our applications, we have omitted a
discussion of the computational aspects; most notably, the division
algorithm and Buchberger's criterion.
\end{remark}

\section{The Gr\"{o}bner Fan}

The purpose of this section is to describe the initial ideals and the
natural adjacency relations among them.  For a commutative polynomial
ring, Mora and Robbiano's Gr\"{o}bner fan \cite{MR} accomplishes this
goal; Assi, Castro-Jim\'{e}nez and Granger~\cite{ACG} have extended
this to the Weyl algebra and we generalize it to the skew polynomial
ring $R$.  Our setting has the advantage that the commutative
polynomial ring, Weyl algebra and homogenized Weyl algebra are all
done at once, shortening the treatment in Saito, Sturmfels and
Takayama~\cite{SST}. 

The next proposition shows that Gr\"{o}bner bases with respect to
vectors $(\uu,\vv)$ generalize those with respect to term orders.

\begin{proposition} 
Let $I$ be an $R$-ideal.  For any term order $\prec$ there exists a
positive vector $(\uu,\vv) \in \PR(R)$ such that $\IN_{\prec}(I) =
\IN_{(\uu,\vv)}(I)$.
\end{proposition}

\begin{proof}
See Proposition~2.1.5 in Saito, Sturmfels and
Takayama~\cite{SST}.
\end{proof}

We prove a key tool in the construction of the Gr\"{o}bner fan.

\begin{proposition} \label{epsilon}
Let $I$ be an $R$-ideal, let $(\uu',\vv')$ belong to $\PR(R)$ and let
$(\uu,\vv)$ be a vector in $\PR(S) = \RR^{m}\times\RR^{n}$.  If
$\varepsilon > 0$ is sufficiently small, then one has
$ 
\IN_{(\uu,\vv)}\big(\IN_{(\uu',\vv')}(I)\big) = \IN_{(\uu'+\varepsilon
\uu,\vv' + \varepsilon \vv)}(I) \, .
$
\end{proposition}

\begin{proof}
Let $\prec$ be any term order and let $\prec_{\varepsilon}$ be the
multiplicative monomial order on $R$ defined as follows:
\[ \begin{array}{rl}
x^{\aa'}y^{\bb'} \prec_{\varepsilon} x^{\aa}y^{\bb} & \text{ if and
only if } \\ & \big(
(\uu'+\varepsilon\uu,\vv'+\varepsilon\vv)\cdot(\aa-\aa',\bb-\bb') > 0
\big) \\ & \text{ or } \left( \begin{array}{l}
(\uu'+\varepsilon\uu,\vv'+\varepsilon\vv)\cdot(\aa-\aa',\bb-\bb') = 0
\\ \text{ and } (\uu',\vv')\cdot(\aa-\aa',\bb-\bb') > 0
\end{array} \right) \\
& \text{ or } \left( \begin{array}{l}
(\uu'+\varepsilon\uu,\vv'+\varepsilon\vv)\cdot(\aa-\aa',\bb-\bb') = 0,
\\ (\uu',\vv')\cdot(\aa-\aa',\bb-\bb') = 0 \\ 
\text{ and } x^{\aa'}y^{\bb'} \prec x^{\aa}y^{\bb} \big)
\end{array} \right).
\end{array} \]
Next, fix a universal Gr\"{o}bner basis $\cU$ for $I$ and choose
$\varepsilon$ small enough so that the following assertions hold:
$(\uu'+\varepsilon \uu,\vv' + \varepsilon \vv) \in \PR(R)$ and for all
elements $g$ in $\cU$, the standard form of $g$ breaks into
four pieces $g(x,y) = g_{0}(x,y) + g_{1}(x,y) + g_{2}(x,y) +
g_{3}(x,y)$ such that $\IN_{\prec_{\varepsilon}}(g) =
g_{0}(\bar{x},\bar{y})$, $\IN_{(\uu'+\varepsilon \uu, \vv'+
\varepsilon \vv)}(g) = g_{0}(\bar{x},\bar{y})+
g_{1}(\bar{x},\bar{y})$, and $\IN_{(\uu', \vv')}(g) =
g_{0}(\bar{x},\bar{y})+ g_{1}(\bar{x},\bar{y})+
g_{2}(\bar{x},\bar{y})$.  In particular, we have
\begin{equation*}
\IN_{(\uu, \vv)}\big(g_{0}(x,y)+g_{1}(x,y)+g_{2}(x,y)\big) =
g_{0}(\bar{x},\bar{y})+ g_{1}(\bar{x},\bar{y}) \, .  \tag{$\dagger$}
\end{equation*}
Since $\cU$ is a Gr\"{o}bner basis with respect to
$\prec_{\varepsilon}$, Proposition~\ref{prec=>(u,v)} provides
two additional Gr\"{o}bner bases in the polynomial ring $S$:
\begin{enumerate}
\item[(i)] The initial forms
$g_{0}(\bar{x},\bar{y})+g_{1}(\bar{x},\bar{y})$ for $g \in \cU$ are a
Gr\"{o}bner basis for the initial ideal $\IN_{(\uu'+\varepsilon
\uu,\vv' + \varepsilon \vv)}(I)$ with respect to
$\prec_{\varepsilon}$.
\item[(ii)] The initial forms $g_{0}(\bar{x},\bar{y}) +
g_{1}(\bar{x},\bar{y}) + g_{2}(\bar{x},\bar{y})$ for $g \in \cU$ are a
Gr\"{o}bner basis for the initial ideal $\IN_{(\uu',\vv')}(I)$ with
respect to $\prec_{\varepsilon}$.
\end{enumerate}
Now, the definition of $\prec_{\varepsilon}$, statement (ii) and
Proposition~\ref{prec=>(u,v)}.2 imply that the polynomials
$g_{0}(\bar{x},\bar{y})+g_{1}(\bar{x},\bar{y})+g_{2}(\bar{x},\bar{y})$
for $g \in \cU$ are a Gr\"{o}bner basis for the ideal
$\IN_{(\uu',\vv')}(I)$ with respect the vector $(\uu,\vv)$.  Moreover,
$(\dagger)$ indicates that the polynomials
$g_{0}(\bar{x},\bar{y})+g_{1}(\bar{x},\bar{y})$ for $g \in \cU$
generate the ideal $\IN_{(\uu,\vv)}\big(\IN_{(\uu',\vv')}(I)\big)$ and
therefore statement (ii) completes the proof.
\end{proof}

We are now in a position to give a description of the equivalence
classes $C_{I}[(\uu,\vv)]$.

\begin{proposition} \label{cone}
Let $I$ be an $R$-ideal, let $\prec$ be a term order and let
$(\uu,\vv)$ belong to $\GR(I)$.  If $\cG$ is the reduced Gr\"{o}bner
basis of $I$ with respect to $\prec_{(\uu,\vv)}$, then one has
\begin{equation*}
C_{I}[(\uu,\vv)] =
\left\{ (\uu',\vv') \in \GR(I) : \IN_{(\uu,\vv)}(g) =
\IN_{(\uu',\vv')}(g) \text{ $\forall$ $g \in \cG$} \right\}
\tag{$\ddagger$}
\end{equation*}
and, hence, each equivalence class $C_{I}[(\uu,\vv)]$ is a relatively
open rational convex polyhedral cone.
\end{proposition}

\begin{proof}
See Proposition~2.3 in Sturmfels~\cite{Stu}.
\end{proof}

We end with the main result of this section.

\begin{theorem}
For $I$ an $R$-ideal, the finite set 
\[
\GF(I) := \left\{ \overline{C_{I}[(\uu,\vv)]} : \text{ for all $(\uu,\vv) \in
\GR(I)$} \right\}
\] 
forms a fan, called the Gr\"{o}bner fan of $I$.
\end{theorem}

\begin{proof}
Given the above lemmas and propositions, the proof is now identical to
Proposition~2.4 in Sturmfels~\cite{Stu}.
\end{proof}

\section{Bounds on the Irreducible Components}

This section contains the proofs of the main results of this paper.
We use the Gr\"{o}bner fan is to show that the Gelfand-Kirillov
dimension of $M$ is independent of the positive vector $(\uu,\vv) \in
\PR(R)$.  We start by stating an elementary lemma from commutative
algebra -- see Lemma~2.2.2 in Saito, Sturmfels and
Takayama~\cite{SST}.

\begin{lemma} \label{dim-drops}
If $J$ is any ideal in $S$ and $(\uu,\vv) \in \RR^{m}\times\RR^{n}$
then one has the following:
\begin{enumerate}
\item[(1)] $\Kdim \IN_{(\uu,\vv)}(J) \leq \Kdim J$.
\item[(2)] If $(\uu,\vv)$ is positive, then
$\Kdim \IN_{(\uu,\vv)}(J) = \Kdim J$. \qed
\end{enumerate}
\end{lemma}

Making use of the Gr\"{o}bner fan, we next study effect that varying
the filtration of $R$ has on the Gelfand-Kirillov dimension of an
$R$-module.

\begin{proposition}  \label{GKdim-indep}
If $M$ is a finitely generated $R$-module. then the Gelfand-Kirillov
dimension of $M$ is independent of the positive vector $(\uu,\vv)
\in\PR(R)$.
\end{proposition}

Recall that $M$ has a good finite dimensional filtration if and only
if the vector $(\uu,\vv)$ is positive.

\begin{proof}
By subsection~8.6.5 in McConnell and Robson~\cite{McR}, we know that
if $M$ has a good finite dimensional filtration, then $\GKdim M =
\GKdim \gr(M)$.  Moreover, the Gelfand-Kirillov dimension and Krull
dimension of $\gr(M)$ are equal by Corollary~\ref{commutative}.  Since
the Krull dimension of a finitely generated module is the Krull
dimension of its support, it suffices to consider ideals.  In
particular, one reduces to proving that, for an $R$-ideal $I$, the
initial ideal $\IN_{(\uu,\vv)}(I)$ is independent of the choice of
positive vector $(\uu,\vv) \in \PR(R)$.

We prove this statement by constructing a homotopy or ``Gr\"{o}bner
walk'' between any two initial ideals.  Let $(\uu_{1},\vv_{1})$ and
$(\uu_{2},\vv_{2})$ be two positive vectors in $\PR(R)$.  We claim
that $\Kdim \IN_{(\uu_{1},\vv_{1})}(I) = \Kdim
\IN_{(\uu_{2},\vv_{2})}(I)$.  Moreover, Proposition~\ref{cone} implies
that each equivalence class $C_{I}[(\uu_{2},\vv_{2})]$ is a convex
cone and, thus, we have $\IN_{(\uu_{2},\vv_{2})}(I) = \IN_{(r \uu', r
\vv')}(I)$ for any positive $r \in \RR$.  Hence, we may, if necessary,
replace $(\uu_{2},\vv_{2})$ by a positive scalar multiple, to
guarantee that $(\uu'-\uu,\vv'-\vv)$ is a positive vector.  It follows
that $(1-r) \cdot (\uu_{1},\vv_{1}) + r \cdot (\uu_{2},\vv_{2})$ is a
positive vector and belongs to $\GR(I)$ for all $r \in [0,1]$.  Define
$J_{r}$ to be the $R$-ideal $\IN_{(1-r) \cdot (\uu_{1},\vv_{1}) + r
\cdot (\uu_{2},\vv_{2})}(I)$.  Since the line segment from
$(\uu_{1},\vv_{1})$ to $(\uu_{2},\vv_{2})$ intersects finitely many
distinct walls of the Gr\"{o}bner fan, there are real numbers $0 =
r_{0} < r_{1} < \cdots < r_{\ell} = 1$ such that the ideal $J_{r}$
remains unchanged as the parameter $r$ ranges inside the open interval
$(r_{j},r_{j+1})$; we denote this ideal by $J_{(r_{j},r_{j+1})}$.  By
Proposition~\ref{epsilon}, we have $\IN_{(\uu_{2},\vv_{2})}(J_{r_{j}})
= J_{(r_{j},r_{j+1})} = \IN_{(\uu_{1},\vv_{1})}(J_{r_{j+1}})$.  By
applying Lemma~\ref{dim-drops} 2, we see that $\Kdim J_{r_{j}} = \Kdim
J_{(r_{j},r_{j+1})} = \Kdim J_{r_{j+1}}$.  Combining these equalities
for $0 \leq j < \ell$ completes the proof.
\end{proof}

We are now ready to prove:

\begin{theorem} \label{dimension_bounds}
Let $p$ be a nonnegative integer and let $M$ be a finitely generated
$R$-module.  If $\Ch_{(\uu,\vv)}(M)$ has an irreducible component of
dimension $p$ for some $(\uu,\vv) \in \PR(R)$, then there is a
submodule $M'$ of $M$ such that $\GKdim M' = p$ when $M'$ is equipped
with a good finite dimensional filtration.
\end{theorem}

\begin{proof}
Fix a positive vector in $\PR(R)$ and let $\tR$ and $\tM$ be
associated Rees ring and module.  Now, suppose $\Ch_{(\uu,\vv)}(M)$
has an irreducible component of dimension $p$.  Since
\[
\Ann_{S}\big(\gr(M)\big) = \gr\big(\Ann_{R}(M)\big) =
\IN_{(\uu,\vv)}\big(\Ann_{R}(M)\big) \, , 
\]
this means that the $S$-ideal $\IN_{(\uu,\vv)}\big(\Ann_{R}(M)\big)$
has a minimal prime $\fp$ of dimension $p$.  By
Proposition~\ref{prec=>(u,v)} and
Proposition~\ref{nontermorder}, we have
\[
\IN_{(\uu,\vv)}\big(\Ann_{R}(M)\big) =
\left. \IN_{(1,\uu,\vv)}\big(\Ann_{\tR}(\tM)\big) \right|_{\bar{x}_{0}=1} \,.
\]   
Because $\bar{x}_{0}$ is a nonzero-divisor on $S[\bar{x}_{0}]$, $\fp$
pulls back to a minimal prime of dimension $p+1$ over
$\IN_{(1,\uu,\vv)}\big(\Ann_{\tR}(\tM)\big)$.  Furthermore, there
exists, by Proposition~\ref{GR=PR}, a positive vector
$(u_{0},\uu',\vv') \in \PR(\tR)$ such that
$
\IN_{(1,\uu,\vv)}\big(\Ann_{\tR}(\tM)\big) =
\IN_{(u_{0},\uu',\vv')}\big(\Ann_{\tR}(\tM)\big) \, .
$
It follows, by Theorem~\ref{GKdim-Kdim}, that $\tM$ has a submodule
$L$ such that $\GKdim L = p+1$.  Finally, Proposition~\ref{pi1dim}
implies $\pi_{0}(L)$ is a submodule of $M$ with Gelfand-Kirillov
dimension $p$.
\end{proof}

Using this theorem, we prove the main results of this paper.

\begin{proof}[Proof of Theorem~\ref{purity}]
If the characteristic variety $\Ch_{(\uu,\vv)}(M)$ has an irreducible
component of dimension strictly less than $p$, then
Theorem~\ref{dimension_bounds} implies that there exists a submodule
$M'$ of $M$ with Gelfand-Kirillov dimension strictly less than $p$.
However, this contradicts our hypothesis.
\end{proof}

We recall Bernstein's inequality which is also called the Weak
Fundamental Theorem of Algebraic Analysis---see Section~1.4 in
Bj\"{o}rk~\cite{bj1} for two distinct proofs.

\begin{theorem}[Bernstein] \label{Bernstein}
Let $k$ be a field of characteristic zero.  If $A_{n}(k)$ has the
standard filtration $(\11,\11)$ and $M$ is a finitely generated
$A_{n}$-module, then one has $\GKdim M \geq n$. \qed
\end{theorem}

\begin{proof}[Proof of Theorem~\ref{main}]
This is immediate from Theorem~\ref{purity} and
Theorem~\ref{Bernstein}.
\end{proof}

\renewcommand{\baselinestretch}{1}


\providecommand{\bysame}{\leavevmode\hbox to3em{\hrulefill}\thinspace}

\end{document}